\documentclass[a4paper]{amsart}
\usepackage{lmodern}

\usepackage[T1]{fontenc}
\usepackage[latin9]{inputenc}
\usepackage{fancyhdr}
\usepackage[english]{babel}
\usepackage{amsmath}
\usepackage{amssymb}

\usepackage[notcite]{showkeys}
\usepackage{hyperref}
\begin{document}

\lhead{\rightmark}

\rhead[\leftmark]{}

\lfoot[\thepage]{}

\cfoot{}

\rfoot{\thepage}

\noindent \begin{center}
\textbf{\Large Solovay's inaccessible over a weak set theory without
choice}
\par\end{center}{\Large \par}

\noindent \begin{center}
{\large Haim Horowitz and Saharon Shelah}%
\footnote{Date: September 9, 2023



Publication 1094 of the second author

Partially supported by European Research Council grant 338821

The second author would like to thank the Israel Science Foundation for partial support of this 

research by grant 1838/19: (ISF) for 2019/10-2023/09) and Rutgers 2018 DMS 1833363: NSF 

DMS Rutgers visitor program (PI S. Thomas) (2018-2022)%
}
\par\end{center}{\large \par}

\noindent \begin{center}
\textbf{Abstract}
\par\end{center}

\noindent \begin{center}
{\small We study the consistency strength of Lebesgue measurability
for $\Sigma^1_3$ sets over Zermelo set theory ($Z$) in a completely choiceless
context. We establish a result analogous to the Solovay-Shelah theorem.}
\par\end{center}{\small \par}

\textbf{\large 0. Introduction}{\large \par}

Our work follows the line of research that was initiated in the celebrated work of Solovay [So] and a later work of the second author [Sh176],  where it was shown that $ZF+DC+$"all sets of reals are Lebesgue measurable" is equiconsistent with $ZFC+$"there exists an inaccessible cardinal". In their works, $DC$ plays an important role. More specifically, the proof in [Sh176] shows that Lebesgue measurability implies $\omega_1^{L[x]}<\omega_1$ for all $x\in \omega^{\omega}$, however, we need a certain amount of choice to conclude that $\omega_1$ is inaccessible in $L$.
\\
\\
One may now ask whether having regularity properties in the complete absence of choice results in no increase in consistency strength. Our main result shows that if we replace $ZF$ with $Z$ (that is, if we remove the Replacement schema), then Lebesgue measurability will still result in an increase in consistency strength in a way analogous to the $ZF+DC$ situation. Namely, we shall prove a version of the following:
\\
\\
\textbf{Theorem} (informal): 1. The following are equiconsistent:
\\
a. $Z+$"all sets of reals are Lebesgue measurable"
\\
b. $ZC+$"there exists an uncountable strong limit cardinal"
\\
\\
2. The following are equiconsistent:
\\
a. $Z+AC_{\aleph_0}+$"all sets of reals are Lebesgue measurable"
\\
b. $ZC+$"there exists a strongly inaccessible cardinal"
\\
\\
Note that the consistency strength of $Z+$"there exists an uncountable strong limit cardinal" is strictly higher than the consistency strength of $Z$, analogously to how $ZFC+$inaccessible has a higher consistency strength than $ZFC$. We also note that the notion of Lebesgue measurability needs a refinement in the complete absence of choice (hence the "informal" in the above statements of the main results), this will require us to work with Borel codes rather than Borel sets, and to consider different notions of measurability that are no longer equivalent in our new setting.
\\
\\
Our first step towards the ultimate goal will be to isolate an ad hoc theory (which we shall denote $Z_*$), which will not include the Replacement schema but will still provide us with the minimal ingredients required in order to imitate the arguments in the Solovay-Shelah theorem. This will be done in Section 1, where we will also isolate some basic relevant measure theoretic notions. We shall then modify and imitate the proofs of the Solovay-Shelah results and obtain their analogs in the context of $Z_*$. This will be done in Sections 2 and 3. As large parts of the work in these sections will consist of modifying the arguments in [So] and [Sh176] to the setting of $Z_*$, we will mainly focus on the parts which are new. For the sake of completeness, we included an appendix with a dense overview of the details of the relevant proof from [Sh176]. Finally, Section 4 will be devoted to transferring our equiconsistency results from $Z_*$ to $Z$. For this purpose, we shall develop in $Z$ a general theory of $L$-like models and show how these can be used to replace $L$ in the equiconsistency proof from Sections 2 and 3.
\\
\\
Finally, we remark that the above-mentioned result opens the door to the following new direction of research:

\textbf{Question}: Given a set theory $T$ that doesn't prove $AC_{\omega}$ and a Suslin
ccc forcing notion $\mathbb Q$,
is $T$ equiconsistent with $T+"$Every set of reals is $\mathbb Q$-measurable$"$?

Our main result provides a negative answer for the case of $T=Z$ and $\mathbb Q=$Random real forcing, but the question remains open for other forcing notions and perhaps over different weak set theories.

\textbf{\large 1. Basic definitions}{\large \par}

Although our main result is formulated for $Z$, we shall first prove our result for a weak theory $Z_*$ (described below), which does not follow from $Z$. We will then show how our arguments translate to a proof over $Z$.

\textbf{Definition 1: }1. $Z^-$ is Zermelo set theory (i.e. $ZF$
without replacement) without choice and without the powerset axiom
(but with the separation scheme).

2. $Z^- + \aleph_n$ is $Z^- +$''the cardinal $\aleph_n$ exists''.

3. $Z_*$ is the theory that consists of the following axioms:

a. $Z^-+\aleph_1 + \aleph_2$.

b. $\mathcal{P}(\mathbb N)$ exists.

c. $L_{\alpha}[z]$ exists for every ordinal $\alpha$ and $z
\in \omega^{\omega}$.

d. $\alpha+\omega$ exists for every ordinal $\alpha$.

\textbf{1.1 Observation ($Z^-$): }Recall that there exists a formula $\phi(x)$
in the language of set theory such that:

a. If $\delta>\omega$ is a limit ordinal, $a\in \omega^{\omega}$
and $b=L_{\delta}[a]$, then $(b,\in) \models \phi(a)$.

b. If $b$ is a transitive set, $a\in \omega^{\omega}$ and $(b,\in) \models \phi(a)$,
then there is a limit ordinal $\delta>\omega$ such that $b=L_{\delta}[a]$.

\textbf{1.2 Convention}: From here until the end of Section 3, our background theory is $Z_*$, so we do not
assume $AC_{\omega}$ and by {}``Borel sets'' we refer only to sets
of reals having a Borel code. In Section 4 we shall translate our results to equiconsistency results over $Z$.

We shall now define several versions of Lebesgue measurability and
the null ideal (note that the different versions are not equivalent
without choice).

\textbf{Definition 2: }1. A set $X\subseteq \mathbb R$ is 1-null
if there exists a Borel set $B$ such that $X\subseteq B$ and $\mu(B)=0$.

2. A set $X\subseteq \mathbb R$ is 2-null if for every $n<\omega$
there exists a Borel set $B_n$ such that $X\subseteq B_n$ and $\mu(B_n)<\frac{1}{n+1}$.

\textbf{Remark: }1-null implies 2-null, and the definitions are equivalent
under $AC_{\omega}$.

\textbf{Definition 3: } A. A set $X\subseteq \mathbb R$ is $i-$measurable
$(i=1,2,3)$ if:

$i=1:$ There exists a Borel set $B$ such that $X\Delta B$ is 1-null.

$i=2:$ There exists a Borel set $B$ such that $X\Delta B$ is 2-null.

$i=3:$ For every $n<\omega$, there exist Borel sets $B_1$ and $B_2$
such that $X\Delta B_1 \subseteq B_2$ and $\mu(B_2) < \frac{1}{n+1}$.

It's easy to see that $i-$measurability implies $j-$measurability
for $i<j$.
\\
\\
B. We define the outer Lebesgue measure $\mu^*$ as usual.
\\
\\
C. For a $3$-measurable set $X$, we let the Lebesgue measure $\mu(X)$ of $X$ be the unique $a\in \mathbb{R}_{\geq 0}$ such that if $n<\omega$, then there are $B_1$ and $B_2$ as above with Borel codes $\eta$ and $\nu$ such that $L_{\omega_1}[\eta, \nu, a] \models "\mu(B_2)< \frac{1}{n+1}$ and $\mu(B_1)-\frac{1}{n+1} \leq a \leq \mu(B_1)+\frac{1}{n+1}"$. This is well-defined when $X$ is $3$-measurable.

\textbf{Observation 4: } a. $X\subseteq \mathbb R$ is 2-null iff $\mu^*(X)=0$. 
\\
\\
b. The Lebesgue measure of a Borel set is well-defined, absolute and does not depend on the Borel code.
\\
\\
c. Analytic statements are absolute between $V$ and $L$, hence $\Sigma^1_2$ statements are upwards absolute from $L$ to $V$.
\\
\\
\textbf{Proof} (of (b)): If $\eta, \nu$ are Borel codes of $B$, work in $L_{\omega_1}[\eta, \nu]$. $\square$ 

\textbf{\large 2. A lower bound on the consistency strength}{\large \par}

\textbf{Claim 5: } Suppose that $V\models Z_*$. 
\\
\\
A.The following version of Fubini's theorem holds:
\\
\\
If $A\subseteq [0,1]$ is not 2-null and $\leq$ is a prewellordering
of $A$ such that every initial segment (i.e. $\{y : y\leq x\}$)
is 2-null, then there exists a set which is not 3-measurable.

B. If in addition $\leq$ is $\Sigma^1_2$, then there exists a non-3-measurable
$\Sigma^1_3$ set. In particular, at least one of the following $\Sigma^1_3$ sets defined in the proof below is not 3-measurable:
$C_{B,a}$, $B_{n,i}$, $A$, $B_0$, $B_1$ and $B_2$.

\textbf{Proof:} For clause (A), suppose that all sets are 3-measurable and we shall
derive a contradiction (for clause (B), it suffices to assume the measurability of the sets mentioned there). For every Borel set $B\subseteq [0,1] \times [0,1]$
and $a\in [0,1]$, define $C_{B,a}:=\{s_1 \in [0,1] : a\leq \mu(\{s_2 : (s_1,s_2) \in B\})\}$,
and similarly, define $C_{a,B}:=\{s_2 \in [0,1] : a\leq \mu(\{s_1 : (s_1,s_2) \in B\})\}$.

Subclaim 1: Let $B\subseteq [0,1] \times [0,1]$ be a Borel set coded
by $r$, and let $a,b,r_1 \in \omega^{\omega}$ such that $r,a,b \in L[r_1]$,
if $L[r_1] \models "\mu(C_{B,a})=b"$ then $V\models "\mu(C_{B,a})=b"$.

Proof: In $L[r_1]$ there is a sequence $(U_n,S_n : n<\omega)$ such
that the sets $U_n \subseteq [0,1]$ are open, the sets $S_n \subseteq [0,1]$
are closed, $S_n \subseteq C_{B,a} \subseteq U_n$ and $L[r_1] \models "\mu(U_n \setminus S_n)<\frac{1}{n}"$.
We shall prove that $\mu(U_n)^V=\mu(U_n)^{L[r_1]}$, $\mu(S_n)^V=\mu(S_n)^{L[r_1]}$
and $V\models "S_n \subseteq C_{B,a} \subseteq U_n"$ for every $n<\omega$ (where $S_n^{L[r_1]}$, etc, are understood using a Borel code).

We shall work in $L[r_1]$ and assume wlog that $S_n \subseteq S_{n+1}$
for every $n<\omega$. Define $R$ as the set of triples $(n,s_1,S)$
such that:

1. $n<\omega$, $s_1 \in S_n \subseteq C_{B,a}$.

2. $S\subseteq [0,1]$ is closed and $\mu(S)=a-\frac{1}{n}$.

3. $s_1 \times S \subseteq B$.

Let $X=\omega \times [0,1]$ and $Y=\{S : S\subseteq [0,1]$ is closed$\}$,
then $X$ and $Y$ are Polish spaces and $R \subseteq X\times Y$
is a $\Pi^1_1$ relation. By $\Pi^1_1$-uniformization, there is a
function $F\subseteq R$ with a $\Pi^1_1$-graph such that $Dom(F)=Dom(R)$.
By absoluteness, the same is true for $(R,F)$ in $V$. Now, if $s_1 \in S_n$
then $s_1 \in S_m$ for every $m>n$ and $\{F(m,s_1) : n\leq m\}$
witnesses that $s_1 \in C_{B,a}$. Therefore, $V\models "sup\{\mu(S_n) : n<\omega\} \leq \mu(C_{B,a})"$.
Similarly we can show that $V\models "\mu(C_{B,a}) \leq inf \{\mu(U_n) : n<\omega\}"$.

Subclaim 2: If $B\subseteq [0,1] \times [0,1]$ is Borel and $a,b\in [0,1]$,
then $\mu(C_{B,a})=b \rightarrow b\leq \frac{\mu(B)}{a}$.

Proof: Let $r_1$ be a real such that $a,b$ and the definition of
$B$ (hence of $C_{B,a}$) are in $L[r_1]$. As the conclusion holds
in $L[r_1]$, it follows from the previous claim that it holds in
$V$ as well.

Subclaim 3: Assume $Z_*$. Fubini's theorem holds for Borel and analytic
sets in the following sense: If $B\subseteq [0,1] \times [0,1]$ is
Borel/analytic and $f_l: [0,1] \rightarrow [0,1]$ $(l=1,2)$ are
defined by $f_l(s_l):=\mu(\{s_{3-l} : (s_1,s_2) \in B\})$, then $\mu(B)=\int f_1(s_1)ds_1=\int f_2(s_2)ds_2$.

Proof: Let $r$ be a real such that the definition of $B$ is in $L[r]$,
and we shall continue the proof as usual in $L[r]$. The only point
that we have to show is that the above integrals are well-defined
and computed in the same way in $L[r]$ and $V$. For every $n>1$
and $i\leq n$, let $B_{n,i}:=\{s_1 : \mu(B_{s_1}^{1}):=\mu(\{s_2 : (s_1,s_2) \in B\}) \in [\frac{i}{n},\frac{i+1}{n}]$.
$(B_{n,i} : i\leq n)$ is a partition of $[0,1]$ for every $n$.
For every $n$, choose a sequence $(S_{n,i},U_{n,i} : i\leq n)$ in
$L[a]$ such that $S_{n,i} \subseteq B_{n,i} \subseteq U_{n,i}$,
$S_{n,i}$ is closed, $U_{n,i}$ is open and $\mu(U_{n,i} \setminus S_{n,i})<\frac{1}{2n}$.
Let $R_1$ be the set of sequences $(n,i,s_1,S)$ such that:

1. $n>1$ and $i\leq n$.

2. $s_1 \in S_{n,i}$.

3. $S\subseteq [0,1]$ is closed and $\frac{1}{n}-\frac{1}{2^n} \leq \mu(S)$.

4. $S\subseteq B_{s_1}^1$.

Let $X:=\{(n,i,s_1) : i\leq n, n>1, s_1 \in S_{n,i}\}$ and $Y$ be
the set of closed subsets of $[0,1]$. As before, by $\Pi^1_1$-uniformization,
there is a $\Pi^1_1$-function $F_1 \subseteq R$ such that for every
$(n,i,s_1)$, if there exists $S$ such that $(n,i,s_1,S) \in R_1$,
then $(n,i,s_1,F_1(n,i,s_1)) \in R_1$. By absoluteness, the same
is true in $V$. Similarly, define $R_2$ as the set of sequences
$(n,i,s_1,U)$ such that:

1. $n>1$ and $i\leq n$.

2. $s_1 \in U_{n,i}$.

3. $U\subseteq [0,1]$ is open and $\mu(U)<\frac{i+1}{n}+\frac{1}{2^n}$.

4. $B_{s_1}^1 \subseteq U$.

As before, there is $\Pi^1_1$ choice function $F_2$ for the relation
$R_2$. $F_1$ and $F_2$ witness that the above integrals are well-defined
and have the same value in $L[r]$ and $V$.

Subclaim 4: If $A \subseteq [0,1]$ and $B=A\times A$, then $\mu^*(B)=\mu^*(A)^2$.

Proof: In one direction, let $a=\mu^*(A)$ and $\epsilon>0$. There
is a Borel set $A^*$ such that $A\subseteq A^* \subseteq [0,1]$
and $\mu(A^*) \leq \mu(A)+\epsilon$. Let $B^*=A^* \times A^*$, then
$\mu^*(B) \leq \mu^*(A^*\times A^*)=\mu(A^*)^2 \leq (a+\epsilon)^2$.
Therefore, $\mu^*(B) \leq \mu^*(A)^2$.

In the other direction, let $a=\mu^*(A)$, $b=\mu^*(B)$ and $\epsilon>0$.
There are Borel sets $A^*$ and $B^*$ such that $A\subseteq A^*$,
$B\subseteq B^*$, $\mu(A^*) \leq \mu^*(A)+\epsilon$ and $\mu(B^*) \leq \mu^*(B)+\epsilon$.
Without loss of generality, $B^* \subseteq A^* \times A^*$. If $s_1 \in A$
then $a\leq \mu(\{s_2 : (s_1,s_2) \in B^*\})$, therefore $a=\mu^*(A) \leq \mu^*(C_{B,a})\leq \mu^*(A^*)<a+\epsilon$.
By Fubini's theorem for Borel sets, it follows that $a^2 \leq \mu^*(B^*)$.
Therefore, $\mu^*(A)^2-\epsilon=a^2-\epsilon \leq \mu^*(B^*)-\epsilon \leq \mu^*(B)$,
so $\mu^*(B)=\mu^*(A)^2$ as required.

We are now ready to complete the proof of claim 5.

Without loss of generality $A\subseteq [0,1]$. We now define the
following sets:

1. $B_0=B=A\times A$

2. $B_1=\{(x,y) \in B : x\leq y\}$

3. $B_2=\{(x,y) \in B : y\leq x\}$

Suppose that each of the sets $A$, $B_0$, $B_1$ and $B_2$ are
3-measurable and we shall derive a contradiction. Choose $\epsilon_1,\epsilon_2 \in (0,1)$
such that $\epsilon_1<\epsilon_2^2$ and $\epsilon_2<\frac{a^2}{6}$
(recall that $A$ is not 2-null by our assumption). As $A$ is 3-measurable,
there are Borel sets $A^*$ and $A^{**}$ such that $A\Delta A^* \subseteq A^{**}$
and $\mu(A^{**})<\epsilon_1$. Similarly, there are Borel sets $B_l^*$
and $B_l^{**}$ $(l=1,2)$ such that $B_l \Delta B_l^* \subseteq B_l^{**}$
and $\mu(B_l^{**})<\epsilon_1$. We shall prove that $\mu(B_l^* \cup B_l^{**})<3\epsilon_2$.
Together we obtain the following:

$a^2=\mu^*(A\times A)=\mu^*(B_1 \cup B_2) \leq \mu^*(B_1^* \cup B_1^{**} \cup B_2^* \cup B_2^{**}) \leq \mu(B_1^* \cup B_1^{**})+ \mu(B_2^* \cup B_2^{**})<3\epsilon_2+3\epsilon_2<a_2$,
a contradiction.

By a previous subclaim, $\mu^*(C_{B_2^{**},\epsilon_2})^V \leq \mu^*(C_{B_2^{**},\epsilon_2})^{L[a]} \leq \frac{\mu^*(B_2^{**})}{\epsilon_2}<\epsilon_2$
(where $a$ is as in the subclaim). Let $C_2:=C_{B_2^{**},\epsilon_2}$
nd let $B_2':=B_2^* \cap ([0,1] \setminus C_2 \times [0,1])$. The
following inequalities hold:

1. $\mu^*(B_2^*) \leq \mu^*(B_2^* \cap (C_2 \times [0,1]))+\mu^*(B_2^* \cap ([0,1] \setminus C_2 \times [0,1]))$

2. $\mu^*(B_2^* \cap (C_2 \times [0,1])) \leq \mu^*(C_2 \times [0,1]) \leq \mu^*(C_2)<\epsilon_2$

Therefore, it suffices to show that $\mu^*(B_2^* \cap ([0,1] \setminus C_2 \times [0,1])) \leq \epsilon_2$.
Given $s_2 \in [0,1] \setminus C_2$, the following holds: $\mu^*(\{s_2 : (s_1,s_2) \in B_2^*\}) \leq \mu^*(\{s_2 : (s_1,s_2) \in B_2\})+ \mu^*(\{s_2 : (s_1,s_2) \in B_2^{**}\}) \leq 0+\epsilon_2$
where the last inequality follows by the choice of $s_1$, the definition
of $B_2$ and the theorem's assumption. By Fubini's theorem, the desired
conclusion follows.

The proof for $l=1$ is similar, where $C_{B_2^{**},\epsilon_2}$
is replaced by $C_{\epsilon_2,B_1^{**}}$ and the rest of the arguments
are changed accordingly. $\square$

\textbf{Theorem 6: }Assume $Z_*$. 

1. If every $\Sigma_3^1$ set of reals is 3-measurable, then $\aleph_1^{L[x]}<\aleph_1$
for every $x\in 2^{\omega}$, hence $\aleph_1$ is a limit cardinal
in $L$.

2. If in addition $AC_{\omega}$ holds, then $\aleph_1$ is inaccessible
in $L$.

\textbf{Proof: }We follow a similar argument as in {[}Sh176{]}. Assume
towards contradiction that $\aleph_1^{L[x_*]}=\aleph_1$ for some
$x_*\in 2^{\omega}$. For every $x\in 2^{\omega}$, let $(\mathbf{B}_{x,i} : i<i(*))$
list all of the Borel null subsets of $2^{\omega}$ (i.e. their Borel
codes, recalling that $"\mu(A)=0"$ is absolute) in $L[x_*,x]$ (we
can do it uniformly in $(x,x_*)$). Denote $=\mathbf{B}_{x,i}^*=\mathbf{B}_{x,i}^V$
and $\mathbf{B}_{x,<i}^*=\underset{j<i}{\cup}\mathbf{B}_{x,j}^*$. Let
$\mathbf{B}_x^*=\underset{i<i(*)}{\cup}\mathbf{B}_{x,i}^*$.

\textbf{Case I: }There exists $x_{**} \in 2^{\omega}$ such that $B_{x_{**}}^*$
is not 2-null.

Work in $V$: Denote $\mathbf{B}=\mathbf{B}_{x_{**}}^*$ and define the
following prewellordering on $\mathbf{B}$: $x\leq y$ iff for every
$i$, $y\in \mathbf{B}_{x_{**},<i}^* \rightarrow x\in \mathbf{B}_{x_{**},<i}^*$.

Cleary, every initial segment of $(\mathbf{B},\leq)$ has the form $\mathbf{B}_{x_{**},<i}^*$,
and hence is 2-null. As $B$ is not 2-null, it follows by claim 5
that there exists a non-3-measurable $\Sigma^1_2$-set, a contradiction.

\textbf{Case II: }$B_x^*$ is 2-null for every $x\in 2^{\omega}$.

We shall first describe the original stages of the proof in {[}Sh176{]},
then we shall describe how to modify the original proof in order to
obtain the desired theorem. The new changes and arguments will be
presented in this section, while the proofs from {[}Sh176{]} will
appear in the appendix.

\textbf{Outline of {[}Sh176{]}:}

We fix a rapidly increacing sequence $(\mu(k) : k<\omega)$ of natural
numbers, say, $\mu(k)=2^{2^{2^{2^{2^{176k}}}}}$.

\textbf{Step I (existence of a poor man generic tree): }Suppose that
$B\subseteq 2^{\omega}$ has measure zero, then there are perfect
trees $T_0,T_1 \subseteq 2^{<\omega}$, functions $m_l: T_l \rightarrow \mathbb{Q}$
and natural numbers $n(k)$ $(k<\omega)$ such that $lim(T_l) \cap B=\emptyset$,
$m_l(\eta)=\mu(lim(T_l) \cap (2^{\omega})^{[\eta \leq]})$ and:

A) 1. $m_0(<>)=\frac{1}{2}$ and for every $\eta \in T_0$, $\mu(lim(T_0) \cap (2^{\omega})^{[\eta \leq]})$
has the form $\frac{k}{4^{lg(\eta)+1}}$ for $0\leq k \leq 4^{lg(\eta)+1}$,
and $k\neq 0$ iff $\eta \in T_0$.

A) 2. $m_1(<>)=\frac{1}{2}$, and for every $\eta \in T_1$, if $lg(\eta) \leq n(k)$
then $\mu(lim(T_1) \cap (2^{\omega})^{[\eta \leq]}) \in \{\frac{l}{4^{n(k)+1}} : 0<l<4^{n(k)+1}\}$.

B) For every $\eta \in 2^{n(k)} \cap T_1$, $2^{n(k)}(1-\frac{1}{\mu(k)})<\mu(lim(T_1) \cap (2^{\omega})^{[\eta \leq]})$.

\textbf{Step II: }Definitions of finite and full systems (see definitions
1-4 in the appendix).

\textbf{Step III: }Showing that the family of finite systems satisfies
ccc (claim 5 in the appendix).

\textbf{Step IV: }Forcing with finite systems over $L[x_*]$ to get
a full system in $L[x_*]$. As the existence of a full system is equivalent
to the existence of a model to a $\mathcal{L}_{\omega_1,\omega}(Q)$
sentence, this is sufficient by absoluteness and Keisler's completeness
theorem.

\textbf{Step V: }We use the full system in order to define two $\Sigma^1_3$
sets of reals (those are the red and the green sets in {[}Sh176{]}),
which will turn out to be non-measurable.

\textbf{Step VI: }Showing that the green and red sets are disjoint,
are not null and have outer measure 1, arriving at a contradiction.

\textbf{Back to the proof of theorem 6:}

We shall describe how each of the above steps should be modified in
order to obtain the proof of our theorem.

\textbf{Step I: Claim: }The claim in step I of $[Sh176]$ holds when
$B$ is a Borel set of measure (say) $<\frac{1}{1000}$. This will
be used in order to show that the red and green sets are not 2-null
(this is step VI).

\textbf{Proof: }Let $r$ be a real that codes $B$. The proof is as
in {[}Sh176{]}, where now we work in $L[r]$. Observe that the tree
$T$ constructed there satisfies $lim(T) \cap A=\emptyset$ where
$A$ is an open set of measure $<\frac{1}{1000}$ containing $B$
(and the construction depends only on $A$).

\textbf{Steps II-III: }No change is needed. 

\textbf{Step IV: }Assuming $Z_*$ we can prove Keisler's completeness
theorem as well as the forcing theorem in $L[r]$ for every $r$ (see
the discussion on forcing over models of $Z_*$ in the end of this
section). Therefore we can repeat the argument in the original Step
IV.

\textbf{Step V: }No change.

\textbf{Step VI: }We shall freely use the notation and definitions
from $[Sh176]$ (see definition 7 and claims 8-11 in the appendix).\textbf{ }

\textbf{Claim A: }The formulas $\phi_{rd}$ and $\phi_{gr}$ are contradictory.

\textbf{Proof: }Suppose that $x$ satisfies both formulas. By definition
7 in the apendix, there is a poor man generic tree over $L[x_*]$
denoted by $T_0^{rd}$ and a poor man generic tree over $L[x_*,T_0^{rd}]$
denoted by $T_1^{rd}$ witnessing $\phi_{rd}(x)$. Repeating the proof
of claim 9 in the appendix, in $L[x_*,x,T_0^{rd},T_1^{rd}]$ there
is a partition $\bar{A}^{rd}=\bar{A}^{rd}(x)=(A_n^{rd} : n<\omega)$
of $\omega_1$ to countably many homogeneously red sets. Similarly,
as $x$ satisfies $\phi_{gr}$, in $L[x_*,x,T_0^{rd},T_1^{rd},T_0^{gr},T_1^{gr}]$
there is a partition $\bar{A}^{gr}=\bar{A}^{gr}(x)=(A_n^{gr} : n<\omega)$
of $\omega_1$ to countably many homogeneously green sets. As $\omega_1$
is regular in $L[x_*,x,T_0^{rd},T_1^{rd},T_0^{gr},T_1^{gr}]$, we
get a contradiction.

\textbf{Claim B: }The formulas $\phi'_{rd}$ and $\phi'_{gr}$ are
contradictory.

\textbf{Proof: }Suppose that $\phi'_{rd}(z) \wedge \phi'_{gr}(z)$,
then for some $x$, $y$ and natual $n^*$ we have $\phi_{rd}(x) \wedge \phi_{gr}(y)$
and $\{n: x(n) \neq y(n)\} \subseteq \{0,...,n^*\}$. Let $\bar{A}^{rd}(x)$
and $\bar{A}^{gr}(y)$ be as in the previous proof, and for every
$n,m<\omega$ let $B_{n,m}=A_n^{rd} \cap A_m^{gr}$. For some $n,m$,
$B_{n,m}$ is infinite. Let $\alpha_k$ be the $k$th element of $B_{n,m}$.
Recalling that $i_1<i_2<i_3 \rightarrow h(i_1,i_2)\neq h(i_2,i_3)$,
then for some $k$ and $j$ we have $h(\alpha_k,\alpha_j)>n^*$. Therefore
$red=x(h(\alpha_k,\alpha_j))=y(h(\alpha_k,\alpha_j))=green$ (recalling
that $\alpha_k,\alpha_j \in A_n^{rd} \cap A_m^{gr}$), which is a
contradiction.

\textbf{Claim C: }$A_{rd}=\{x : \phi_{rd}(x)\}$ and $A_{gr}=\{x : \phi_{gr}(X)\}$
are not of measure zero.

\textbf{Proof: }This is the same argument as in claim 10 in the appendix,
the only difference is that instead of taking a $G_{\delta}$ set
of measure zero covering $A_{rd}$, we take for every $0<\epsilon$
a Borel set of measure $<\epsilon$ covering $A_{rd}$. By the modified
construction of the poor man generic tree, we continue as in the original
proof. 

\textbf{Claim D: }$A_{rd}$ is not 3-measurable.

\textbf{Proof: }As in $[Sh176]$ (claim 11 in the appendix). $\square$

\textbf{\large 3. An upper bound on consistency strength (following Levy)}{\large \par}

Historical remark: While Solovay's proof used the Levy collapse of an inaccessible
cardinal (which results in a model of $DC$), our proof follows an
older argument of Levy that used the collapse of a limit uncountable
cardinal.

\textbf{Theorem 7: }$A\rightarrow B$ where:

A) 1. $V\models Z_*C$.

2. $V=L$.

3. $\lambda$ is a limit cardinal $>|P(\mathbb{N})|$ such that $\mu<\lambda \rightarrow 2^{\mu}<\lambda$.

4. $\mathbb{P}=\Pi \{\mathbb{P}_{\mu,n} : \mu<\lambda, n<\omega\}$
is a finite support product such that $\mathbb{P}_{\mu,n}=Col(\omega,\mu)$. 

5. $G\subseteq \mathbb{P}$ is generic, $\eta_{\mu,n}=\underset{\sim}{\eta_{\mu,n}}[G]: \omega \rightarrow \mu$
is the generic of $\mathbb{P}_{\mu,n}$.

6. In $V[G]$ we define $V_1=V[\{\eta_{\mu,n} : \mu<\lambda,n<\omega\}]$,
i.e. the class of sets in $V[G]$ hereditarily definable from parameters
in $V$ and a finite number of members of $\{\eta_{\mu,n} : \mu<\lambda,n<\omega\}$.

B) 1. $V_1 \models Z_*$.

2. $V_1 \models \aleph_1=\lambda$.

3. If $\lambda$ is singular in $V$ then then $V_1 \models cf(\lambda)=\aleph_0$.

4. If $\lambda$ is regular in $V$ then $V_1 \models cf(\lambda)=\aleph_1$.

5. The following claim holds in $V_1$: If (a)+(b)+(c) hold then (d)
holds where:

a. $\mathbb{Q}$ is a defnition of a forcing notion (with elements
which are either reals or belong to $H(\aleph_1)$) with parameters
in $V_1$ satisfying c.c.c., such that $\mathbb{Q}$ is absolute enough
in the following sense: There is $\bar{t}_*=((\mu_i,n_i) : i<n(*))$
such that $\mathbb{Q}$ is definable using $\bar{\eta_{t_*}}=\{\eta_{\mu_i,n_i} : i<n(*)\}$
and parameters from $V$, and if $\bar{t}=((\mu_l,n_l) : l<n)$ then
$\mathbb{Q}^{V[\bar{\eta}_{t_*t}]} \lessdot \mathbb{Q}^{V_1}$.

b.1. $\underset{\sim}{\eta}$ is a $\mathbb{Q}$-name of a real, i.e.
a sequence of $\aleph_0$ antichains given in $V[\bar{\eta}_{\bar{t}_*}]$.

b.2. The generic set can be constructed from $\underset{\sim}{\eta}$
in a Borel way.

c. The ideal $I=I_{(\mathbb{Q},\underset{\sim}{\eta}),\aleph_0}$
(see {[}HwSh1067{]}) satisfies: $\bar{t}_* \leq \bar{t}_1 \leq \bar{t}_2 \rightarrow I^{V[\bar{\eta}_{\bar{t}_1}]}=P(P(\mathbb{N}))^{V[\bar{\eta}_{\bar{t}_1}]} \cap I^{V[\bar{\eta}_{\bar{t}_2}]}$.

Remark: Note that $P(\mathbb{N})^{V_1}=\cup \{P(\mathbb{N})^{V[\bar{\eta}_{\bar{t}}]} : \bar t$
has the form $((\mu_i,n_i) : i<n)\}$.

d. Every $X\subseteq \omega^{\omega}$ equals a Borel set modulo $I$.

We shall first outline Solovay's original proof from {[}So{]}, then
we shall describe how to smilarly prove the above theorem.

\textbf{An outline of Solovay's proof (for random real forcing)}

\textbf{Step I: }Let $G\subseteq Coll(\omega,<\kappa)$ be generic
where $\kappa$ is inaccessible and let $x\in V[G] \cap Ord^{\omega}$,
then there exists a generic $H\subseteq Coll(\omega,<\kappa)$ such
that $V[G]=V[x][H]$.

\textbf{Step II: }For every formula $\phi$ there is a formula $\phi^*$
such that for every $x\in V[G] \cap Ord^{\omega}$, $V[G]\models \phi(x)$
iff $V[x]\models \phi^*(x)$.

\textbf{Step III: }In $V[G]$, $\omega^{\omega} \cap V[a]$ is countable
for every $a\in Ord^{\omega}$.

\textbf{Step IV: }For every $a\in \omega^{\omega}$, $\{x\in \omega^{\omega} : x$
is not $(\mathbb{Q},\underset{\sim}{\eta})$-generic over $V[a]\} \in I$,
where $\mathbb Q$ is random real forcing and $\underset{\sim}{\eta}$
is the name for the generic.

\textbf{Step V: }Given a maximal antichain $J\in V[a]$ of closed
sets deciding $\phi^*(a,\underset{\sim}{\eta})$ (where $\underset{\sim}{\eta}$
is the name for the random real), we define the desired Borel set
as union of members of $J$ forcing $\phi^*(a,\underset{\sim}{\eta})$.

\textbf{Proof of theorem 7: }Suppose that $A\subseteq \omega^{\omega}$
is definable using $\bar{\eta}_{\bar{t}}$ for $\bar{t}=((\mu_i,n_i) : i<n)$.
As before, we shall indicate how to modify Solovay's original proof
for our purpose.

\textbf{Step I: }Our aim is to prove a result similar to Step I above,
where the real parameter belongs to $V_1$.\textbf{ }Suppose that
$a\in V_1$ is a real (so $a=\underset{\sim}{a}[G]$ for some $\mathbb{P}$-name
$\underset{\sim}{a}$), then $a$ is definable by a formula $\phi$
from a finite number of $\eta_{\mu,n}$'s, say $\{ \eta_{\mu_i,n_i} : i<i(*)\}$.
In order to prove that $a\in V[\{ \eta_{\mu_i,n_i} : i<i(*)\}]$,
it's enough to show that:

Claim 7.1: If $p\in \mathbb{P}$ and $p\Vdash \underset{\sim}{a}(n)=k$
then $p\restriction \underset{i<i(*)}{\Pi}\mathbb{P}_{\mu_i,n_i} \Vdash \underset{\sim}{a}(n)=k$.

Proof: Suppose towards contradiction that $p\restriction \underset{i<i(*)}{\Pi}\mathbb{P}_{\mu_i,n_i} \leq q$
forces a different value for $\underset{\sim}{a}(n)$. Let $\pi$
be an automorphism of $\mathbb{P}$ over $\underset{i<i(*)}{\Pi}\mathbb{P}_{\mu_i,n_i}$
such that $\pi(p)$ is compatible with $q$ (just switch the relevant
coordinates), then $\pi(p) \Vdash \underset{\sim}{a}(n)=k$, a contradiction. 

In order to complete this step, we shall prove the following claim:

Claim 7.2: If $\mathbb{Q} \lessdot \underset{i<i(*)}{\Pi}\mathbb{P}_{\mu_i,n_i}$
then there is an isomorphism of $RO(\mathbb{P})$ onto $RO(\mathbb{Q} \times \mathbb{P})$
that is the identity over $RO(\mathbb{Q})$.

Proof: Let $\kappa>\mu>\aleph_1+max\{\mu_i : i<i(*)\}$. As we assume
that $V=L$ (so in particular we have GCH), the usual proof works.

Conclusion 7.3: If $G\subseteq \mathbb{P}$ is generic over $V$ and $a\in (\omega^{\omega})^{V_1}$,
then there is a generic $H\subseteq \mathbb{P}$ such that $V[G]=V[a][H]$.

Proof: By the above claims, $a\in V[\{\eta_{\mu_i,n_i} : i<i(*)\}]$
for an appropriate finite set of $\eta_{\mu_i,n_i}$'s. Let $B_a$
be the the complete subalgebra generated by $a$, then by the previous
claim $B_a \times \mathbb{P}$ is isomorphic to $\mathbb{P}$ (over
$B_a$) and the claim follows.

\textbf{Steps II: }Same as in Solovay's proof.

\textbf{Step III:} Suppose that $G' \subseteq \underset{i<i(*)}{\Pi}\mathbb{P}_{\mu_i,n_i}$
is generic. We shall use the fact that the ideal $I$ is generated
by sets which are disjoint to some $B_N$ such that $N\subseteq H(\aleph_1)^{L[G']}$
where: 

1. $N$ is transitive and $||N||=\aleph_0$.

2. $B_N=\{\underset{\sim}{\eta}[H] : H$ is $\mathbb{Q}^{L[G']} \cap N$-generic
over $N\}$.

Work in $V_1$: Let $X$ be the set of $\nu \in \omega^{\omega}$
(in $V_1$) such that $\nu$ is not $(N,\mathbb{Q},\underset{\sim}{\eta})$-generic
where $N=(H^{V[G']}(\aleph_1),\in)$. As $N$ is countable in $V_1$
(recall that $\lambda$ is strong limit) and $\nu \in (\omega^{\omega})^{V_1}$
is generic over $N$ iff it's generic over $V[G']$, it follows by
the definition of $I$ that $X\in I$.

\textbf{Step IV: }Suppose that $A\subseteq \omega^{\omega}$ is definable
by $\phi(\bar{\eta}_{\bar t},x)$ where $\bar{\eta}_{\bar t} \in L[G']$
and $G' \subseteq \underset{i<n}{\Pi}\mathbb{P}_{\mu_i,n_i}$ is the
generic set obtained by the restriction of $G$ to $\underset{i<n}{\Pi}\mathbb{P}_{\mu_i,n_i}$.
Let $\{p_n : n<\omega\}\subseteq \mathbb{Q}^{L[G']}$ be a maximal
antichain and $(\underset{\sim}{\mathbf{t}_n} : n<\omega)$ a sequence
of names of truth values such that $p_n \Vdash \underset{\sim}{\eta} \in A$
iff $\underset{\sim}{\mathbf{t}_n}=true$ (such sequences exist by step
II). Let $\bar{p}=(\bar{p}^i : i<\omega)$ enumerate all maximal antichains
in $H(\aleph_1)^{L[G']}$ (so each $\bar{p}^i$ is of the form $\bar{p}^i=(p_n^i : n<\omega)$). 

By our assumption, given a generic real $\eta$ we can define the
set $G_{\eta}$ in a Borel way such that: 

$(*)$ $G_{\eta}$ is generic over $H(\aleph_1)^{L[G']}$ and $\underset{\sim}{\eta}[G_{\eta}]=\eta$.

Now let $B:=\{ \eta : G_{\eta}$ is well-defined, satisfies $(*)$
above and for some $n$, $p_n \in G_{\eta} \wedge \underset{\sim}{\mathbf{t}_n}[G_{\eta}]=true \}$.
Denote by $B_n$ the set of $\eta\in B$ such that $"\eta \in B"$
is witnessed by $n$.

$B$ is Borel by our assumptions on the forcing. Therefore it's enough
to prove that $A=B$ $mod$ $I$.

Let $\eta \in \omega^{\omega}$ (in $V_1$), by step III it's enough
to show that if $\eta$ is generic over $H(\aleph_1)^{L[G']}$ (and
hence $\eta=\underset{\sim}{\eta}[G_{\eta}]$ for $G_{\eta}$ as in
$(*)$ above) then $\eta \in A$ iff $\eta \in B$. Indeed, if $\eta \in A$
(and $\eta=\underset{\sim}{\eta}[G_{\eta}]$ where $G_{\eta}$ is
as in $(*)$), by the definition of $\{p_n : n<\omega\}$ and $(\underset{\sim}{\mathbf{t}_n} : n<\omega)$,
there is some $p_n \in G_{\eta}$ such that $\underset{\sim}{\mathbf{t}_n}[G_{\eta}]=true$,
therefore $\eta \in B_n \subseteq B$. Similarly, if $\eta \in B_n$
for some $n$ such that $\underset{\sim}{\mathbf{t}_n}[G_{\eta}]=true$,
then by the definitions of $\{p_n : n<\omega\}$ and $(\underset{\sim}{\mathbf{t}_n} : n<\omega)$,
$\eta \in A$. $\square$

\textbf{Conclusion 8: }A) The following theories are equiconsistent
for $i\in \{1,2,3\}$:

1. $Z_*C+$''there is a limit cardinal $>\aleph_0$''.

2. $Z_*C+$''there is a strong limit cardinal$>\aleph_0$''.

3$(i)$. $Z_*+$''every $\Sigma_3^1$ set of reals is $i$-measurable''.

4$(i)$. $Z_*+$''every set of reals is $i$-measurable.

B) The following theories are equiconsistent for $i\in \{1,2,3\}$:

1. $Z_*C+$''there is a regular limit cardinal$<\aleph_0$''.

2. $Z_*C+$''there is strongly inaccessible cardinal''.

3$(i)$. $Z_*+DC+$''every $\Sigma_3^1$ set of reals is $i$-measurable''.

4$(i)$. $Z_*+DC+$''every set of reals is $i$-measurable''.

5$(i)$. $Z_*+AC_{\aleph_0}+$''every set of reals is $i-$measurable''.
\\
\\
\textbf{Proof}: By putting together all the above. $\square$

\textbf{\large A remark on forcing over models of $Z_*$}{\large \par}

In order to guarantee that the generic extensions in our proofs satisfy
$Z_*$, we work in the context of models of $Z_*$ of the form $L$
or $L[r]$ for some real $r$. In this context, we work with classes
$W$ of the following form: There is a formula $\phi$ with parameters
that defines the class, and there is a limit ordinal $\nu<\omega^2$
such that $\phi$ defines $W\cap L_{\alpha}[r]$ in $L_{\alpha+\nu}[r]$
when $\alpha$ is a limit ordinal (recall that for every ordinal $\alpha$,
the ordinal $\alpha+\omega n$ exists).

Now, for a set forcing $\mathbb P$ in a model of the above form,
we define the class of $\mathbb P-$names as above. Therefore, for
every limit ordinal $\alpha$ we define the intersection of $L_{\alpha}[r]$
with the class of names. For the names that we defined, we can prove
the forcing theorem as usual and show that $Z_*$ holds in the generic
extension. In addition, note that when we force over $L[r]$, as $L[r]$
has a well-ordering $<_{L[r]}$ definable from $r$, we can use it
to get a well-ordering of the generic extension, hence a model of
$Z_*C$.

\textbf{4. Translating the proofs from $Z_*$ to $Z$}

Our goal in this section will be to prove a version of Corollary 8 for $Z$. Note that $Z$ implies all axioms in Definition 1, except of (c) and (d), and so we our first step will be to provide an adequate substitute for them that can be established in $Z$. As we can't prove in $Z$ that $\alpha+\omega$ exists for every ordinal $\alpha$, we shall avoid using the von Neumann definition of an ordinal and instead refer by an "ordinal" to the order type of a well-ordered set as defined below. We may avoid using proper classes with the help of the following definition:
\\
\\
\textbf{Definition 9}: a. Let $n<\omega$. A $V_{\omega+n}$-ordinal is the order type of a well-ordered set whose set of elements is contained in $V_{\omega+n}$. Pedantically, it's a $E_{\omega + n}$-equivalence class where $E_{\omega + n}$ is the equivalence relation consisting of all pairs $((A_1, <_1),(A_2, <_2))$ such that, for $l=1, 2$, $A_l \subseteq V_{\omega + n}$, $(A_l, <_l)$ is a well ordering, $(A_1, <_1) \cong (A_2, <_2)$ and $(A_l, <_l)$ is not isomorphic to a well order $(B, <)$ where $B\subseteq V_{\omega +m}$ for some $m<n$.
\\
\\
b. We say that $\alpha$ is a $V_{< \omega + \omega}$-ordinal or a $*$-ordinal if it's a $V_{\omega + n}$-ordinal for some $n<\omega$. $*$-ordinals will be denoted by $\alpha, \beta, \gamma$, etc. We say that $\alpha$ is a $V_{<\omega+n}$-ordinal if it's a $V_{\omega + m}$-ordinal for some $m<n$.
\\
\\
c. The natural ordering $\leq_*$ on $*$-ordinals will be defined by $\alpha \leq_* \beta$ iff every $(A, <)$ in $\alpha$ is isomorphic to an initial segment of some $(B, <')$ in $\beta$. 
\\
\\
d. Let $Ord_{\omega+\omega}$ denote the class of $V_{<\omega+\omega}$-ordinals.
\\
\\
e. For $n<\omega$, let $Ord_{\omega+n}$ denote the set of $\alpha$ such that for some $m\leq n$, $\alpha$ is a $V_{\omega+m}$ ordinal. 
\\
\\
f. For $\alpha \in Ord_{\omega+ \omega}$, let $set(\alpha):=\{ \beta : \beta <_* \alpha\}$.
\\
\\
\textbf{Observation 10}: a. If $(A,<)$ is a well-order and $A\subseteq V_{\omega+n}$, then $(A, <) \in V_{\omega +n+6}$, $(A,<)/E_{\omega+n} \in V_{\omega+n+7}$ and $set((A,<)/E_{\omega+n}) \in V_{\omega+n+8}$.
\\
\\
b. If $\alpha$ is a $V_{\omega+n}$-ordinal, then $\alpha \in V_{\omega+n+7}$ and $set(\alpha) \in V_{\omega+n+8}$ (and so both are sets).
\\
\\
\textbf{Proof}: (b) follows directly from (a) and the definition of $V_{\omega+n}$-ordinals. As for (a), if $x,y \in A \subseteq V_{\omega+n}$, then $(x, y) \in V_{\omega+n+2}$, and so the ordering $<$ is in $V_{\omega+n+3}$ and in a similar fashion we establish the rest of the claim. $\square$
\\
\\
\textbf{Definition 11}: Let $\alpha$ be a $*$-ordinal and $A\subseteq \alpha$.
\\
\\
a. Fix a vocabulary $\tau=\{ \in, P\}$ where $\in$ and $P$ (which we intend to interpret as $A$) are a binary and unary predicates, respectively. We let $\Phi$ be the set of all first-order formulas $\phi=\phi(x, \bar{x}_u)$ in the vocabulary $\tau$ such that $u\subseteq \omega$ is finite and $x$, $\bar{x}_u=(x_n : n\in u)$ are the free variables of $\phi$.
\\
\\
b. Given a $*$-ordinal $\beta$, let $\Psi_{\beta}=\{ (u, \phi(x, \bar{x}_u), f)$: $u$ and $\phi(x, \bar{x}_u)$ are as in (a), $f$ is a function from $u$ to $\beta \} \cup set(\beta)$ (note that this is a disjoint union).
\\
\\
\textbf{Remark 11A}: Pedantically, it should be noted at this point that while $Z$ doesn't prove the existence of $V_{\omega}=\cup_{n<\omega} V_n$, this is inconsequential for our purposes, as any infinite set as $V_{\omega}$ will do. Therefore, we may define $V_{\omega}$ to be $\omega$.
\\
\\
\textbf{Observation 12}: Assume $Z$.
\\
\\
a. For every $n<\omega$, the set $V_{\omega +n}$ exists, and $\{(x, n) : x\in V_{\omega + n}\}$ is a definable class.
\\
\\
b. For every $n<\omega$, the set $Ord_{\omega+n}$ exists, as well as $<_* \restriction Ord_{\omega + n}$ and the function $\alpha \mapsto set(\alpha)$ for $\alpha \in Ord_{\omega+n}$. $Ord_{\omega+\omega}$, $<_*$ and the function $\alpha \mapsto set(\alpha)$ for $\alpha \in Ord_{\omega + \omega}$ are definable classes.
\\
\\
c. $\Psi_{\alpha}$ exists for every $*$-ordinal $\alpha$, and $\{( \beta, \Psi_{\beta}) : \beta$ is a $*$-ordinal$\}$ is a definable class.
\\
\\
d. If $\beta \in Ord_{\omega+n}$ and $\alpha \leq_* \beta$, then $\alpha \in Ord_{\omega+n}$.
\\
\\
e. There is a fixed $k<\omega$ (say, $k=20$) such that, for $n<\omega$, there is a unique element $\alpha_n^* \in Ord_{\omega+n+k}$ such that $set(\alpha_n^*)=Ord_{\omega+n}$. It follows that $\alpha_n^* \notin Ord_{\omega+n}$ and $\alpha_n^* \leq_* \alpha_{n+1}^*$.
\\
\\
f. $\alpha_n^*$ is a $*$-cardinal, that is, $\beta <_* \alpha_n^* \rightarrow |set(\beta)|< |set(\alpha_n^*)|$.
\\
\\
\textbf{Proof}: Clauses (a)-(d) are straightforward. For example, (c) follows by the existence of products and power sets, followed by an application of Separation. For clause (e), let $(A_1, <_1)=(Ord_{\omega+n}, <_*)$. We saw that $A_1 \subseteq V_{\omega+n+7}$. Let $m<\omega$ be minimal such that there are $A_2 \subseteq V_{\omega+m}$ and a well-order $<_2$ of $A_2$ such that $(A_2, <_2) \cong (A_1, <_1)$, then $\alpha_n^*:=(A_2, <_2)/E_{\omega+m}$ is as required. Finally, we prove clause (f). Suppose that $\beta <_* \alpha_n^*$, then $set(\beta) \subseteq set(\alpha_n^*)$, so $|set(\beta)| \leq |set(\alpha_n^*)|$. As $set(\alpha_n^*)=Ord_{\omega+n}$, it follows that $\beta \in Ord_{\omega+n}$. Therefore, there is $m\leq n$, $A\subseteq V_{\omega+m}$  and a well ordering $<_1$ of $A$ such that $\beta=(A,<_1)/E_{\omega+m}$. Let $k<\omega$, $B\subseteq V_{\omega+n+k}$ and $<_2$ be a well-order of $B$ such that $(B, <_2)/E_{\omega+n+k}=\alpha_n^*$. So $(B, <_2) \cong (set(\alpha_n^*), <_* \restriction set (\alpha_n^*))$ and $(A, <_1) \cong (set(\beta), <_* \restriction set(\beta))$. Suppose towards contradiction that $|set(\beta)|=|set(\alpha_n^*)|$, so there is a bijection $h: B\rightarrow A$ and so we get the existence of the ordering $<':=\{ (h(a), h(b)) : a,b \in B, a<_2 b\}$. Therefore, $(A, <') \cong (set(\alpha_n^*), <_* \restriction set (\alpha_n^*))$, and as $A\subseteq V_{\omega+m}$, we get $\alpha_n^* \in Ord_{\omega+n}$, contradicting the previous clause. It follows that $|set(\beta)|<|set(\alpha_n^*)|$. $\square$
\\
\\
\textbf{Definition 13}: Given a $*$-ordinal $\beta$ and $A\subseteq set(\beta)$, let $\mathcal{L}_{\beta, A}$ be the set of all objects $\mathbf m$ of the form $\mathbf m=(X, O, R, P, L, E)=(X_{\mathbf m}, O_{\mathbf m}, R_{\mathbf m}, P_{\mathbf m}, L_{\mathbf m}, E_{\mathbf m})$ such that:
\\
\\
a. $X=\Psi_{\beta}$.
\\
\\
b. $E$ is an equivalence relation on $X$ such that $\gamma /E=\{ \gamma \}$ for every $\gamma <_* \beta$.
\\
\\
c. $O\subseteq X$ is $set(\beta)$ (this is intended to be the set of ordinals in our model), and we may identify it with $set(\beta)/E$.
\\
\\
d. $R$ is a binary relation on $X/E$ such that, for $\gamma_1, \gamma_2 \in set(\beta)$, $\{\gamma_!\}R\{\gamma_2\}$ iff $\gamma_1 <_* \gamma_2$.
\\
\\
e. $P=A/E \subseteq X/E$. 
\\
\\
f. $(X/E, R)$ satisfies a large enough finite fragment of $Z^-$.
\\
\\
g. $L$ is a binary relation on $X/E$ (the intention is that $(x/E, y/E) \in L$ should correspond to $x/E \in L_{y/E}$).
\\
\\
h. $(x/E, y/E) \in L \rightarrow y/E \in O$.
\\
\\
i. $(x_1 /E, x_2/E)\in R \wedge (x_2 /E, y/E) \in L \rightarrow (x_1 /E, y/E) \in L$.
\\
\\
j. $(\{x/E : (x/E, y/E) \in L  \} : y/E \in O/E)$ is increasing continuous.
\\
\\
k. If $y/E$ is the first element in $O$, then $\{x/E : (x/E, y/E) \in L \}=\emptyset$.  
\\
\\
l. If $\{ \delta \}=z/E \in O/E$ is the successor of $\{ \gamma \}=y/E$ and $(s/E, z/E) \in L$, then $S:=\{x/E : (x/E, s/E) \in R\}$ is a subset of $\{x/E : (x/E , y/E) \in L\}$ of the following form: Let $N$ be the $\tau$-model with universe $\{ t/E : (t/E, y/E) \in L \}$, $\in^N = R \restriction N \times N$ and $P^N=A$, then there is some $(u, \phi(x, \bar{x}_u), f) \in \Psi_{\gamma}$ such that $S$ is definable over $N$ using $\phi(x,...,f(i)/E,...)_{i\in u}$.
\\
\\
m. For every $y, z$ and $(u, \phi(x, \bar{x}_u), f)$ as above, there is a corresponding $s$ as above.
\\
\\
\textbf{Claim 14}: Assume $Z$.
\\
a. For every $\beta$ and $A\subseteq set(\beta)$, $\mathcal{L}_{\beta, A}$ is non-empty.
\\
\\
b. Moreover, $\mathcal{L}_{\beta, A}$ is a singleton, denoted $\mathcal{L}_{\beta, A}=\{ \mathbf{m}_{\beta, A}\}$.
\\
\\
c. $\mathcal{L}_{\leq \alpha, A}:= \cup \{ \mathcal{L}_{\beta, A} : \beta \leq \alpha\}$ and $\{ (\beta, \mathbf{m}) : \beta \leq \alpha, \mathbf m \in \mathcal{L}_{\beta, A}\}$ exist.
\\
\\
\textbf{Proof}: By Power Set and Separation, with the existence and uniqueness claims following a similar line as in the case of the standard construction of the constructible universe. $\square$
\\
\\
\textbf{Definition 15}: We shall now define our $L$-like models in the context of $Z$. Note that if $\alpha <_* \beta$, $A \subseteq set(\beta)$ and $\mathbf{m} \in \mathcal{L}_{\beta, A}$, then $\mathbf m \restriction \alpha$ is naturally defined and is the unique member of $\mathcal{L}_{\alpha, A \cap \alpha}$. For $\gamma<_* \beta$, let $L^{\dagger}_{\gamma}[A]$ be the model $N$ from Definition 13(l). 
\\
\\
\textbf{Observation 16}: In Definition 15, $L_{\gamma}^{\dagger}[A]$ doesn't depend on $\beta$. $\square$
\\
\\
We shall assume WLOG that the natural numbers of $V$ and $L^{\dagger}_{\gamma}[A]$ are the standard natural numbers.
\\
\\
\textbf{Definition 17}: Given $\alpha \in Ord_{\omega+ \omega}$ and $A\subseteq set(\alpha)$, let $L^{\dagger}[A]= \cup \{ L^{\dagger}_{\gamma}[A] : \gamma$ is a $*$-ordinal$\}$.
\\
\\
\textbf{Claim 18}: Assume $Z$. 
\\
\\
A) a. Let $m<n$ and $A\subseteq set(\alpha_m^*)$. If $B\subseteq L_{\gamma}^{\dagger}[A]$ for some $\gamma<_* \alpha_m^*$ and $B\in L_{\alpha_n^*}^{\dagger}[A]$, then $B\in L_{\alpha_m^*}^{\dagger}[A]$. 
\\
\\
b. There are infinitely many cardinals in $L^{\dagger}[A]$.
\\
\\
B) Let $\alpha_*$ be a $*$-ordinal and $A\subseteq set(\alpha_*)$, then:
\\
\\
a. $L^{\dagger}[A] \models Z$.
\\
\\
b. $L^{\dagger}[A] \models "\alpha+ \omega$ exists for every ordinal $\alpha$".
\\
\\
c. $L^{\dagger}[A] \models Z_*$.
\\
\\
\textbf{Proof}: Part A:
\\
a. By a similar argument as in Goedel's proof of GCH in $L$, using Observation 12(f).
\\
\\
b. Follows from Observation 12.
\\
\\
Part B:
\\
a. The only nontrivial axioms are the Power Set and Separation axioms. For the Power Set axiom, let $n_* < \omega$ such that $\alpha_* \in Ord_{\omega+n_*}$. Let $X \in L^{\dagger}[A]$, then there is some $n>n_*$ such that $X\in L_{\alpha_n^*}^{\dagger}[A]$. For a $*$-ordinal $\beta$ with $\alpha_n^*<_* \beta$, let $P_{\beta}(X)=\{ y\in L^{\dagger}_{\beta}[A] : L^{\dagger}_{\beta}[A] \models "y\subseteq X"\}$, then $P_{\beta}(X) \in L_{\beta+1}^{\dagger}[A]$, so it suffices to show that $P_{\beta}(X)$ stabilizes for large enough $\beta$, which now follows from clause (A)(a). 
\\
It remains to establish the Separation schema in $L^{\dagger}[A]$. Let $X\in L^{\dagger}[A]$, $\bar a\in L^{\dagger}[A] ^{<\omega}$ and $\phi(x, \bar y)$ be a formula such that $\phi(x, \bar a)$ defines a subset of $X$. Let $\delta \in Ord_{\omega+\omega}$ such that $range(\bar a) \subseteq L_{\delta}^{\dagger}[A]$ and $X\in L_{\delta}^{\dagger}[A]$. As $L^{\dagger}[A]$ has unboundedly many cardinals, there is $\lambda$ such that $L^{\dagger}[A] \models "\delta \leq_* \lambda$ and $\lambda$ is a cardinal$"$. Let $F_0,...,F_{n-1}$ list the Skolem functions for the existential subformulas of $\phi$ and for the formula asserting that a model $M$ has the form $L^{\dagger}_{\alpha}[A]$. Let $Y$ be the closure of $L^{\dagger}_{\lambda}[A]$ inside $L^{\dagger}[A]$ under $F_0,...,F_{n-1}$ (this is a subclass of $L^{\dagger}[A]$ and $L^{\dagger}[A] \restriction Y \models "V=L^{\dagger}[A]"$). Let $W$ be the class of all triples $(\epsilon, \zeta, f)$ such that, in $L^{\dagger}[A]$:
\\
a. $\epsilon <_* \lambda^+$ is a $*$-ordinal with $\lambda \leq_* \epsilon$.
\\
b. $\zeta$ is a $*$-ordinal and $\epsilon \leq_* \zeta$.
\\
c. $f$ is an isomorphism from $L^{\dagger}_{\zeta}[A] \cap Y$ to $L^{\dagger}_{\epsilon}[A]$ (and so is the identity on $L_{\lambda}^{\dagger}[A]$). 
\\
Note that if $(\epsilon_i, \zeta_i, f_i) \in W$ ($i=1,2$), then $f_1 \subseteq f_2$ or $f_2 \subseteq f_1$. Let $\epsilon(*)=sup \{ \epsilon \leq_* \lambda^+ : (\epsilon, \zeta, f)\in W \}$ and let $f_*=\cup \{ f : (\epsilon, \zeta, f) \in W$ and $\epsilon<\epsilon(*) \}$.
\\
Now $X\in L_{\epsilon(*)}^{\dagger}[A]$, $f_* \restriction X=id$ and $f_*(\bar a)= \bar a$. By the choice of $f_*$, its range is $L^{\dagger}_{\epsilon(*)}[A]$. The subclass of $X$ definable in $L_{\epsilon(*)}^{\dagger}[A]$ by $\phi(x, \bar a)=\phi(x, f_*(\bar a))$ is in $L^{\dagger}_{\epsilon(*)+1}[A]$ hence in $L^{\dagger}[A]$. It remains to show that for $b\in L^{\dagger}_{\lambda}[A]$, $L^{\dagger}_{\epsilon(*)}[A] \models \phi(b, \bar a)$ iff $L^{\dagger}[A] \models \phi(b, \bar a)$. As $Y$ is closed under the relevant Skolem functions, $L^{\dagger}[A] \models \phi(b, \bar a)$ iff $L^{\dagger}[A] \restriction Y \models \phi(b, \bar a)$. As $f_*$ is an isomorphism from $Dom(f_*)$ to $L^{\dagger}_{\epsilon(*)}[A]$ which is the identity on $b$ and $\bar a$, it will suffice to show that $Dom(f_*)=Y$. Suppose towards contradiction that there is some $y\in Y \setminus Dom(f_*)$. $y\in L^{\dagger}_{\xi}[A]$ for some $*$-ordinal $\xi$. Choose a $<_*$-minimal $*$-ordinal $\xi$ for which there is such a pair $(y, \xi)$, so necessarily $\xi=\zeta+1$ is a successor $*$-ordinal $<_*$-greater than $\lambda$. By Separation in $V$, $Y\cap L^{\dagger}_{\xi}[A]$ is a set. We shall obtain a contradiction by constructing an isomorphism from $Y\cap L^{\dagger}_{\xi}[A]$ to $L^{\dagger}_{\epsilon}[A]$ for some $\epsilon \leq_* \epsilon(*)$. For this purpose, it will suffice to prove the following general subclaim:
\\
\textbf{Subclaim}: Assume $V\models Z$. For $l=1,2$, (A)($l$) implies (B) where:
\\
\\
A($l$). a. Let $\tau=\{R, S, P\}$ where $R$ and $S$ are two-place predicates and $P$ is a unary predicate.
\\
\\
b. For $l=1,2$, $M_l=(|M_l|, R^{M_l}, S^{M_l}, P^{M_l})=(A_l, R_l, S_l, P_l)$ are $\tau$-models satisfying:
\\
\\
- $(P_l, R_l)=(P_l, R_l \restriction P_l)$ is an infinite linear well-order.
\\
- For $ a\in A_l$,  letting $A_{l,a}^1:=\{b\in A_l : bR_la \}$, we have $A_{l, a_1}^1=A_{l, a_2}^1 \rightarrow a_1=a_2$. 
\\
- $S_l \subseteq \{ (a,b): a\in A_l, b\in P_l \}$.
\\
- For $b\in P_l$, let $A_{l, b}^2:= \{ a: (a, b) \in S_l \}$. 
\\
- For $a_1, a_2 \in P_l$, we have $a_1 S_l a_2 \rightarrow A_{l, a_1}^2 \subseteq A_{l, a_2}$.
\\
- If $a\in P_l$ is limit in $(P_l, R_l)$, then $A_{l, a}^2=\cup \{ A_{l,b}^2 : b\in P_l, bR_l a\}$.
\\
- If $a\in P_l$ is first in $(P_l, R_l)$, then $A_{l, a}^2=\emptyset$.  
\\
- If $c\in A_l$, then either $c\in A_{l, a}^2$ for some $c\in P_l$ or $A_{l, c}^1 \subseteq A_{l,a}^2$ and $a$ is the last member of $(P_l, R_l)$.
\\
- If $a,b\in P_l$ and $b$ is the successor of $a$ in $(P_l, R_l)$, then:
\\
In the case of (A)($1$), $\{A_{l, c}^1 : c\in A_{l,b}^2\}$ is the family of first-order definable subsets of $A_{l, a}^2$. In the case of (A)($2$), $A_{l,b}^2$ is the union of $A_{l, a}^2$ and the set of all subsets of $A_{l,a}^2$ obtained by a Goedel operation.   
\\
\\
B. One of the following holds:
\\
\\
a. $M_1$ and $M_2$ are isomorphic.
\\
\\
b. There is $a_2 \in P_2$ such that $M_1$ is isomorphic $M_2 \restriction A_{2, a_2}^2$.
\\
\\
c. There is $a_1 \in P_1$ such that $M_1 \restriction A_{1, a_1}^2$ is isomorphic to $M_2$.
\\
\\
In order to prove the subclaim, we shall assume for simplicity that $(P_l, R_l)$ don't have a last element. Let $W$ be the set of all triples $(\epsilon, \zeta, f)$ such that $\epsilon \in A_1$, $\zeta \in A_2$ and $f$ is an isomorphism from $M_1 \restriction A_{1, \epsilon}^2$ to $M_2 \restriction A_{2, \zeta}^2$. As before, if $(\epsilon_i, \zeta_i, f_i) \in W$ $(i=1,2)$, then $f_1 \subseteq f_2$ or $f_2 \subseteq f_1$. Now let $B_1=\{ \epsilon \in P_1: (\epsilon, \zeta, f) \in W$ for some $\zeta$ and $f \}$ and $B_2=\{ \zeta \in P_2: (\epsilon, \zeta, f) \in W$ for some $\epsilon$ and $f \}$. If $B_1=P_1$ and $B_2=P_2$, then $M_1$ is isomorphic to $M_2$. If $B_1=P_1$ and $B_2 \neq P_2$, then letting $a_2 \in P_2 \setminus B_2$ be minimal, we have that $M_1$ is isomorphic to $M_2 \restriction A_{2, a_2}^2$ (and similarly in the case where $B_1 \neq P_1$ and $B_2=P_2$). Finally, we observe that it's impossible to have $B_1 \neq P_1$ and $B_2 \neq P_2$, as in this case we choose minimal $a_l \in P_l \setminus B_l$ and use the last property in clause (A)($l$)(b) to extend the isomorphism thus obtaining a contradiction. This concludes the proof of the subclaim.
\\
We can now apply the (proof of) the subclaim with $M_1 \restriction A_{1, \epsilon}^2=L_{\epsilon}^{\dagger}[A]$ and $M_2 \restriction A_{2, \zeta}^2=Y\cap L_{\zeta}^{\dagger}[A]$ (both of which are sets). We note that the models $M_2 \restriction A_{2, \zeta}^2=Y\cap L_{\zeta}^{\dagger}[A]$ satisfy the assumptions of the subclaim. The only nontrivial part is the last assumption, for which we may use the case (A)($l$) for $l=2$ and further stratify the $L_{\alpha}^{\dagger}[A]$-hierarchy, letting $L_{\alpha,0}^{\dagger}[A]=L_{\alpha}^{\dagger}[A]$, $L_{\alpha, n+1}^{\dagger}[A]$ be the union of $L_{\alpha, n}^{\dagger}[A]$ and the subsets obtained by a Goedel operation and $L_{\alpha+1}^{\dagger}[A]=\cup_{n<\omega} L_{\alpha,n}^{\dagger}[A]$. This completes the proof of Separation in $L^{\dagger}[A]$ and thus the proof of clause (B)(a).
\\
\\
b. Follows from the fact that the $*$-ordinals are closed under addition with $\omega$.
\\
\\
c. Follows from (b). $\square$
\\
\\
In order to prove our final result, we shall generalize the above results and definitions from $*$-ordinals to more general well-ordered classes and sets.
\\
\\
\textbf{Definition and Observation 19}: Let $S$ and $<$ be classes such that $(S, <)$ is a well-order and such that $S_{<a}:=\{ b\in S : b<a\}$ is a set for every $a\in S$. Given a bounded set $A\subseteq S$, we can define the model $L^{\ddagger}[A,S]$ analogously to the way we defined the models $L^{\dagger}[B]$ (with $B$ a bounded set of $*$-ordinals). That is, for every $a\in S$ we can define the set $\Psi_a$ as we did in Definition 11 and repeat the construction in Definition 13 with $S$ replacing $Ord_{\omega+\omega}$, $(\Psi_a : a\in S)$ replacing $(\Psi_{\alpha} : \alpha \in Ord_{\omega+\omega})$ and $<$ replacing $<_*$. We observe that the analogs of Claim 14 and Observation 16 hold in the same way for this generalization.
\\
\\
Observe that $L^{\dagger}[A]$ coincides with $L^{\ddagger}[A, S]$ for $(S, <)=(Ord_{\omega+\omega}, <_*)$. We also note that the property of $Ord_{\omega+\omega}$ that was used to prove Claim 18 is the fact that the cardinals in $L^{\dagger}[A]$ are unbounded. This motivates the following definition:
\\
\\
\textbf{Definition 20}: A triple $(S, <, A)$ as in Definitin 19 is called good if $L^{\ddagger}[A,S]$ has unboundedly many cardinals.
\\
\\
\textbf{Observation 21}: If $(S,<,A)$ is good, then $L^{\ddagger}[A,S]$ satisfies all of the properties of the models $L^{\dagger}[B]$ established above, in particular, $L^{\ddagger}[A,S] \models ZC$. $\square$
\\
\\
\textbf{Claim 22}: a. $(Ord_{\omega+\omega}, <_*, A)$ is good for any bounded set $A\subseteq Ord_{\omega+\omega}$.
\\
\\
b. If $V\models "\lambda$ is an uncountable limit cardinal$"$ and $A\subseteq \lambda$ is a bounded set, then $(\lambda, <, A)$ is good.
\\
\\
c. If $V\models "\lambda$ is an uncountable limit cardinal$"$, then there is a good triple $(S, <, \emptyset)$ such that $\lambda \in S$ and $L^{\ddagger}[\emptyset, S] \models "\lambda$ is an uncountable limit cardinal$"$.
\\
\\
\textbf{Proof}: (a) has already been established and (b) is obvious. As for clause (c), let $V_{[0]}:=\lambda$, and for $n<\omega$, let $V_{n+1}:=P(V_{[n]})$. Now define $V_{[n]}$-ordinals replacing $V_{\omega+n}$ by $V_{[n]}$ in the definition of $V_{\omega+n}$-ordinals. Denote the resulting class by $(Ord_{[\omega]},<_{**})$. We can now repeat the exact same arguments as in the case of $*$-ordinals to establish that $(Ord_{[\omega]},<_{**},\emptyset)$ is good. Note that pedantically $\lambda$ in $Ord_{[\omega]}$ is identified with its equivalence class. It can now be checked that $\lambda$ is an uncountable limit cardinal in $L^{\ddagger}[Ord_{[\omega]}, <_{**}, \emptyset]$. $\square$
\\
\\
Finally, before proving our equiconsistency results, we observe the existence of an uncountable limit cardinal is a large cardinal axiom over $Z$:
\\
\\
\textbf{Observation 23}: The consistency strength of $Z+"$there exists an uncountable limit cardinal$"$ is higher than the consistency strength of $Z$.
\\
\\
\textbf{Proof}: Suppose that $V\models "\lambda$ is a limit uncountable cardinal$"$. By Observation 22(b), $(S,<, A)=(\lambda, <, \emptyset)$ is good, so $L^{\ddagger}[A,S] \models ZC$. Now note that $L^{\ddagger}[A,S]=L_{\lambda}^{\ddagger}[A,S]$ is a set. $\square$
\\
\\
Combining the above results with the previous sections, we thus arrive to our final conclusion:
\\
\\
\textbf{Conclusion 24}: A) The following theories are equiconsistent
for $i\in \{1,2,3\}$:

1. $ZC+$''there is a limit cardinal $>\aleph_0$''.

2. $ZC+$''there is a strong limit cardinal $>\aleph_0$''.

3$(i)$. $Z+$''every $\Sigma_3^1$ set of reals is $i$-measurable''.

4$(i)$. $Z+$''every set of reals is $i$-measurable.

B) The following theories are equiconsistent for $i\in \{1,2,3\}$:

1. $ZC+$''there is a regular limit cardinal $>\aleph_0$''.

2. $ZC+$''there is strongly inaccessible cardinal''.

3$(i)$. $Z+DC+$''every $\Sigma_3^1$ set of reals is $i$-measurable''.

4$(i)$. $Z+DC+$''every set of reals is $i$-measurable''.

5$(i)$. $Z+AC_{\aleph_0}+$''every set of reals is $i-$measurable''.
\\
\\
C) The following theories are equiconsistent
for $i\in \{1,2,3\}$:
\\
\\
1. $ZF$
\\
\\
2. $ZF+$"every $\Sigma^1_3$ set of reals is $i$-measurable".
\\
\\
3. $ZF+$"every set of reals is $i$-measurable".
\\
\\
\textbf{Proof}: Throughout the proof, we shall freely use the fact that, assuming $Z$, $L^{\dagger}$ and $L^{\dagger}[r]$ $(r\in 2^{\omega})$ are models of $Z$ and $Z_*$.
\\
\\
\textbf{(A)}: We split the proof to two cases. In each case, we shall show that if $V$ is a model of one clause, then we can construct a model of any of the other clauses. 
\\
\textbf{Case I}: There is a good triple $(S, <, A)$ such that $L^{\ddagger}[A,S]$ has a limit uncountable cardinal.
\\
In this case, $L^{\ddagger}[A,S]$ is a model of clauses (1) and (2) (as $L^{\ddagger}[A,S] \models GCH$). In order to get the consistency of clauses (3) and (4), we repeat the proof from Section (3) over $L^{\dagger}[A,S]$.
\\
\textbf{Case II}: For every good triple $(S,<,A)$, $L^{\ddagger}[A,S]$ has no limit uncountable cardinal.
\\
It will suffice to show in this case that $V$ is not a model of any of the clauses. We shall first show that $V$ is not a model of clause $3(i)$ (hence also not of clause $4(i)$). By the assumption of the claim, $L^{\dagger}=L^{\ddagger}(Ord_{\omega+\omega},<_*, \emptyset)$ has no limit uncountable cardinal. Let $\Gamma:= \{ \alpha \in Ord_{\omega+\omega}:$ every $(A,<) \in \alpha$ is countable$\}$, then $\Gamma \subseteq Ord_{\omega}=set(\alpha_0^*)$. Furthermore, $\Gamma$ is a $<_*$-downward closed subset of $set(\alpha_0^*)$. Therefore, there is $\gamma_1^* \leq \alpha_0^*$ such that $\Gamma=set(\gamma_1^*)$ (and so $\omega^{L^{\dagger}} <_* \gamma_1^*$). As $L^{\dagger}$ has no limit uncountable cardinal, it follows that there is some $\gamma_0^* <_* \gamma_1^*$ such that $L^{\dagger} \models "\gamma_0^*$ is a cardinal and $\gamma_1^*=(\gamma_0^*)^+ "$. Now let $r\subseteq \omega$ code a well-ordering of $\omega$ order type $\gamma_0^*$ (recalling that $\omega^V=\omega^{L^{\dagger}})$, then $L^{\dagger}[r] \models "\gamma_1^*=\aleph_1"$. Finally, note that there is a natural translation of reals from $V$ to reals from $L^{\dagger}$: Let $S:=\{(\eta, \nu): \eta \in (2^{\omega})^V, \nu \in (2^{\omega})^{L^{\dagger}[r]}$ and for every $n<\omega$,  $V\models \eta(n)=1$ iff $L^{\dagger}[r] \models \nu(n)=1 \}$. $S$ exists in $V$ by the usual arguments (i.e., as $V$ and $L^{\dagger}[r]$ satisfy the Power Set axiom, so we can form the relevant product and use Separation), and it gives rise to a one-to-one function from $(2^{\omega})^{L^{\dagger}[r]}$ to $(2^{\omega})^V$. Furthermore, suppose that $L^{\dagger}[r] \models "B$ is a definition of a Borel subset of $2^{\omega}"$, so there is a corresponding well-founded subtree $T \subseteq \omega^{<\omega}$ and a function $h$ with domain $T$ such that $h(\eta)$ is a clopen subset of $2^{\omega}$ if $\eta \in T$ is minimal, and $h(\eta) \in \{ \cup, \cap\}$ otherwise. Then we can naturally interpret $B$ as a Borel set in $V$. We are now in the same setting as in Theorem 6 from Section (2), and we can repeat the proof there to obtain a non-measurable $\Sigma^1_3$ set. It follows that $V$ is not a model of clauses $3(i)$ and $4(i)$. Finally, the fact that clauses $(1)$ and $(2)$ don't hold in $V$ follows from Claim 22(c).
\\
\\
\textbf{(B)}: Clause (B) is similar, splitting the cases according to whether or not there is a good triple $(S,<,A)$ such that $L^{\ddagger}[A,S]$ has an inaccessible cardinal. If there is such a triple, we repeat the arguments of Case I with the use of Section (3) replaced by Solovay's original argument over $L^{\ddagger}[A,S]$ so we can also obtain $DC$. Suppose now that there is no good triple $(S,<,A)$ such that $L^{\ddagger}[S,<,A]$ has an inaccessible cardinal. Note that if $V$ satisfies clause (1) or (2), then for the good triple $(S,<,\emptyset)$ from Claim 22(c) we have $L^{\ddagger}[\emptyset, S] \models "\lambda$ is an inaccessible cardinal$"$, contradicting the assumption of our case. Therefore, $V$ doesn't satisfy clauses (1) or (2). For the remaining clauses, it suffices to show that $V$ is not a model of $Z+AC_{\aleph_0}+\Sigma^1_3$-Lebesgue measurability. If $V \models \neg AC_{\aleph_0}$, then we are done. So assume that $V\models Z+AC_{\aleph_0}$. By our assumption, $L^{\dagger}$ has no inaccessible cardinals. By $AC_{\aleph_0}$, $\gamma_1^*$ is regular, and therefore must be a successor in $L^{\dagger}.$ We now proceed as in clause (A).
\\
\\
\textbf{(C)}: Finally, the non-trivial part of (C) follows by the same argument as in Section 3. $\square$
\\
\\
\\
\\

\textbf{\large Appendix: Can you take Solovay's inaccessible away?
({[}Sh176{]})}{\large \par}

We now copy the definitions, theorems and proofs from {[}Sh176{]}
that are relevant for understanding the above proofs. 

\textbf{The following definitions are presented as step II on page 6.}

\textbf{Definition 1. }1. Let $N_n$ be the set of pairs $(t,m)$
such that:

a. $\emptyset \neq t\subseteq 2^{\leq n}$ is closed under initial
segments, and for every $\eta \in t\cap 2^{<n}$, for some $l$, $\eta \hat <l> \in t$.

b. $m: t\rightarrow \mathbb{Q}$ is a function such that $m(<>)=\frac{1}{2}$,
$4^{lg(\eta)+1}m(\eta) \in \mathbb{N} \cap [1,4^{lg(\eta)+1}2^{-lg(\eta)})$,
and for $\eta \in t\cap 2^{<n}$, $m(\eta)=\Sigma\{m(\eta \hat <l>) : \eta \hat <l> \in t\}$.

2. Let $N=\underset{n<\omega}{\cup}N_n$, we call $n$ the height
of $(t,m)$ for $(t,m)\in N_n$ and denote it by $ht(t,m)$. If $t'=t \cap 2^{\leq n}$,
$m'=m \restriction t'$, we let $(t',m')=(t,m) \restriction n$. There
is a natural tree structure on $N$ defined by $(t_0,m_0)\leq (t_1,m_1)$
if $(t_0,m_0)=(t_1,m_1) \restriction ht(t_0,m_0)$.

3. A closed tree $T\subseteq 2^{<\omega}$ satisfies $(t,m)$ if $T\cap 2^{\leq ht(t,m)}=t$
and for every $\eta$, $\mu(lim(T) \cap (2^{\omega})_{[\eta]})=m(\eta)$.

\textbf{Definition 2. }1. $M_k$ is the set of pairs $(t,m)$ such
that for some $n=ht(t,m)$ we have:

a. $\emptyset \neq t\subseteq 2^{\leq n}$ is closed under initial
segments, and for $\eta \in t\cap 2^{<n}$ there is $l\in \{0,1\}$
such that $\eta \hat <l> \in t$.

b. $m: t\rightarrow \mathbb{Q} \cap (0,1)$ is a function such that
$m(<>)=\frac{1}{2}$, and for $\eta \in t\cap 2^{<n}$, $m(\eta)=\Sigma\{m(\eta \hat <l>) : \eta \hat <l> \in t\}$.

c. We define $r_l=lev_l(t,m)$ by induction on $l$: $r_0=0$, $r_{i+1}$
is the first $r>r_i$ such that $r\leq n$, for every $\eta \in 2^{\leq r} \cap t$,
$4^{r+1}m(\eta) \in \mathbb{N}$, and for every $\eta \in 2^r \cap t$,
$m(\eta)>2^{-r}(1-\frac{1}{\mu(l+1)})$.

Now we demand that $r_k$ is well defined and equals $n$.

2. Let $M_{k,n}=\{(t,m)\in M_k : ht(t,m)=n\}$, $M_{k,<n}=\underset{l<n}{\cup}M_{k,l}$,
$M_{*,<n}=\underset{k<\omega}{\cup}M_{k,<n}$, $M=\underset{k<\omega}{\cup}M_k$. 

3. For $(t,m) \in M_k$, let $rk(t,m)=k$.

4. We define the order on $M$ as we did for $N$.

\textbf{Definition 3: }A finite (full) system $S$ consists of the
following:

A. The common part: A finite subset $W\subseteq \omega_1$ (the set
$W=\omega_1$) and a number $n(1)<\omega$ ($n(1)=\omega$) and a
function $h:[W]^2 \rightarrow n(1)$ such that if $i_1<i_2<i_3$ belong
to $W$, then $h(i_1,i_2)\neq h(i_2,i_3)$.

B. The red part:

a. For every $(t,m) \in M_{*,\leq n(1)}$ there is a natural number
$\lambda(t,m)$, and for every $(t_1,m_1) \in N_{\lambda(t,m)}$ there
is a member $\rho(t_1,m_1,t,m) \in t\cap 2^{ht(t,m)}$.

b. Let $\{\eta_l : l<\omega\}$ be a fixed enumeration of $2^{<\omega}$
such that $lg(\eta_l) \leq l$. For every $(t,m) \in M_{k,\leq n(1)}$,
$l<k$, $j<k$ and $\xi \in W$, there is a finite set $A_{l,j}^{(t,m),\xi} \subseteq 2^{\leq \lambda(t,m)}$
such that $\underset{\nu \in A_{l,j}^{(t,m),\xi}}{\Sigma}\frac{1}{2^{lg(\nu)}}<\frac{1}{2^{l+j}}$.

c. For every $(t,m) \in M_{k,\leq n(1)}$, $\xi \in W$ and $(t(0),m(0)) \in N_{\lambda(t,m)}$
there is a function $f_{(t(0),m(0))}^{(t,m),\xi}: \{\eta_l : l<k\} \times k \rightarrow \omega$.

d. Monotonicity for (a): If $(t_0,m_0)<(t_1,m_1)$ (both in $M_{*,\leq n(1)}$),
then $\lambda(t_0,m_0)<\lambda(t_1,m_1)$. Moreover, if $(t^0,m^0)<(t^1,m^1) \in N_{\lambda(t_1,m_1)}$,
then $\rho(t^0,m^0,t_0,m_0)<\rho(t_1,m_1,t_1,m_1)$.

e. Monotonicity for (b): If $(t^0,m^0)<(t^1,m^1)$ (both in $M_{*,\leq n(1)}$)
and $A_{l,j}^{(t^0,m^0),\xi}$ is defined, then $A_{l,j}^{(t^0,m^0),\xi}=A_{l,j}^{(t^1,m^1),\xi}$.
Also $f_{(t_0,m_0)}^{(t^0,m^0),\xi} \subseteq f_{(t_1,m_1)}^{(t_1,m_1),\xi}$
if $(t_0,m_0)<(t_1,m_1) \in N_{\lambda(t^1,m^1)}$.

f. The homogeneity consistency condition: If $(t,m) \in M_{k,\leq n(1)}$,
$\xi<\zeta \in W$, $h(\xi,\zeta)<ht(t,m)$, $(t_1,m_1) \in N_{\lambda(t,m)}$
and $\rho=\rho(t_1,m_1,t,m)$, then:

1. $\rho(h(\xi,\zeta))=0(=red)$

or

2. For every $l,j<k$, $j\neq 0$ such that $f_{(t_1,m_1)}^{(t,m),\zeta}(\eta_l,j)=f_{(t_1,m_1)}^{(t,m),\xi}(\eta_l,j)$
there is no perfect tree $T\subseteq 2^{<\omega}$ which satisfies
$(t_1,m_1)$ and $t_1^{\eta_l \leq}$ is disjoint to $\underset{\alpha<k}{\cup}A_{\alpha,j}^{(t,m),\xi}$
and to $\underset{\alpha<}{\cup}A_{\alpha,j}^{(t,m),\zeta}$.

C. The green part: It is defined similarly, only in (f)(1) we replace
$0(=red)$ by $1(=green)$.

\textbf{Definition 4. }The order between finite systems is defined
naturally (for a given $(t,m)$, $\lambda(t,m)$, $A_{l,j}^{(t,m),\xi}$,
$f_{(t(0),m(0))}^{(t,m),\xi}$ remain fixed, $W$ and $n(1)$ might
become larger).

\textbf{The following claim corresponds to step III on page 6.}

\textbf{Claim 5: }The family of finite systems satisfies the countable
chain condition.

\textbf{Proof: }Let $(S(\gamma) : \gamma<\omega_1)$ be a sequence
of $\omega_1$ conditions. By a delta-system argument, we may assume
that for $S(0)$ and $S(1)$ we have: $n:=n(1)^{S(0)}=n(1)^{S(1)}$,
$\lambda^{S(0)}=\lambda^{S(1)}$, $\rho^{S(0)}=\rho^{S(1)}$ and there
is a bijection $g: W^{S(0)} \rightarrow W^{S(1)}$ such that $g$
is the identity on $W^{S(0)} \cap W^{S(1)}$ and $g$ maps $S(0)$
onto $S(1)$ in a natural way.

We shall define a common upper bound $S$. We let $W^S:W^{S(0)} \cup W^{S(1)}$,
$n(1)^S=n+1$. The function $h^s$ is defined as follows: By the above
claim, we may assume that $h^{S(0)}$ agrees with $h^{S(1)}$ on $W^{S(0)} \cap W^{S(1)}$.
$h^s$ will extend $h^{S(0)} \cup h^{S(1)}$ as follows: If $\xi<\zeta \in W^S$
and $\xi \in W^{S(l)} \iff \zeta \notin W^{S(l)}$ $(l=0,1)$, then
$h^S(\xi,\zeta)=n$. For each $(t,m) \in M_{*,\leq n}$ we let $\lambda(t,m)$,
$\rho(-,-,t,m)$ be as in $S(0)$ and $S(1)$, and for $\xi \in W^{S(l)}$,
$A_{l,j}^{(t,m),\xi}$ and $f_{(t_1,m_1)}^{(t,m),\xi}$ are defined
as in $S(l)$.

We shall now define the above information for $(t,m) \in M_{*,leq n+1} \setminus M_{*,\leq n}$.
So let $(t,m) \in M_{k+1,\leq n+1} \setminus M_{*,\leq n}$, hence
$ht(t,m)=n+1$. Clearly there is a unique $(t(0),m(0))<(t,m)$, $(t(0),m(0)) \in M_{*,\leq n}$
$(M_{k,\leq n})$. WLOG we shall concentrate on the red part. Define
$\lambda(t,m)=\lambda((t(0),m(0)))+|W^S|+(2k+1)$. For every $j\leq k$
define an independent family $(A_{k,j}^{(t,m),\xi} : \xi \in W^S)$
of subsets of $\{\nu : lg(\nu)=\lambda(t,m)\}$ such that $\frac{|A_{k,j}^{(t,m),\xi}|}{2^{\lambda(t,m)}}=\frac{1}{2^{k+j+1}}$.

Define $f_{(t_1,m_1)}^{(t,m),\xi}(\eta_l,j)$ for $(t_1,m_1) \in N_{\lambda(t,m)}$,
$j,l<k+1$ as follows:

1. If $j,l<k$, $\xi \in W^{S(l)}$, let $f_{(t_1,m_1)}^{(t,m),\xi}(\eta_l,j)=f_{(t_1,m_1) \restriction \lambda(t_0,m_0)}^{(t_0,m_0),\xi}(\eta_l,j)$.

2. If $l=k$ or $j=k$, we think of $f_{(t_1,m_1)}^{(t,m),\xi}(\eta_l,j)$
as a function of $\xi$, and we shall define it arbitrarily as an
injective function to $\omega$ (recalling that $W^S$ is finite).

Defining $\rho(t_1,m_1,t,m)$ for $(t_1,m_1)\in N_{\lambda(t,m)}$:

Let $(t_0,m_0):=(t_1,m_1) \restriction \lambda(t(0),m(0))$ (by monotonicity,
$\lambda(t(0),m(0))<\lambda(t,m)$) and $\rho_2=\rho(t_0,m_0,t(0),m(0)) \in t(0)$,
so $lg(\rho_2)=ht(t(0),m(0))$. We shall find a proper extension $\rho \in t$
of $\rho_2$ that will satisfy definition $3(f)$. We shall consider
the cases where $3(f)(2)$ fails, in each such case we need to guarantee
that $\rho(h(\xi,\zeta))=0$. Now $lg(\rho_2)=ht((t(0),m(0)))$. Recall
that $(t(0),m(0))\in M_{k,\leq n} \subseteq M_k$, therefore, $r_k$
in definition 2(c) exists and equals $ht(t(0),m(0))=\lambda(t(0),m(0))$.
By the definition of $M_k$, $m(0)(\rho_2)>2^{-ht(t(0),m(0))}(1-\frac{1}{\mu(k)})$.
By 2(b), $m(0)(\rho_2)=\underset{\rho_2 \leq \nu \in t\cap 2^{n+1}}\sum m(0)(\nu)$,
now suppose that $|\{ \nu \in t\cap 2^{n+1}\} : \rho_2 \leq \nu| \leq 2^{(n+1)-lg(\rho_2)}(1-\frac{1}{\mu(k)})$,
then $m(0)(\rho_2) \leq 2^{-ht(t(0),m(0))}(1-\frac{1}{\mu(k)})$ as
$m(0)(\nu) \leq 2^{-(n+1)}$ for every $\rho_2 \leq \nu \in t\cap 2^{n+1}$,
which is a contradiction. Therefore, $|\{ \nu \in t\cap 2^{n+1}\} : \rho_2 \leq \nu| > 2^{(n+1)-lg(\rho_2)}(1-\frac{1}{\mu(k)})$.
Therefore, if $3(f)(2)$ fails for less than $log(\mu(k))$ quadruples,
then we can find $\rho$ that satisfies the demands in $3(f)$ (suppose
not, then for some $c<log(\mu(k))$, there are $c$ coordinates above
$lg(\rho_2)$ such that no extension of $\rho_2$ in $t$ of length
$n+1$ has $0$ in those coordinates. There are $\frac{2^{n+1-lg(\rho_2)}}{2^c}$
sequences with 0 in those coordinates, therefore, the number of extensions
in $t$ is at most $2^{n+1-lg(\rho_2)}(1-\frac{1}{2^c})$ which is
a contradiction).

For a given pair $(l,j)$ we want to count the number of $\xi \in W^S$
such that $t_1^{[\eta_l \leq]}$ is disjoint to $\underset{\alpha<k+1}{\cup}A_{\alpha,j}^{(t,m),\xi}$.
Our goal is to show that $3(f)(2)$ fails for $<log(\mu(k))$ choices
of $(l,j,\xi,\zeta)$. Now recall that $lg(\eta_l) \leq l$, and by
the definition of $N$, $\frac{1}{4^{l+1}} \leq \frac{1}{4^{lg(\eta_l)+1}} \leq m_1(\eta_l)$.
As before, $\frac{|t_1^{[\eta_l \leq]}|}{2^{\lambda(t(0),m(0))}-lg(\eta_l)}>\frac{1}{4^{l+1}}$.
Recall that $|\{\nu \in A_{k,j}^{(t,m),\xi} : \eta_l \leq \nu\}|=\frac{|A_{k,j}^{(t,m),\xi}|}{2^{lg(\eta_l)}}$,
therefore, if $x$ is the number of sets $A_{k,j}^{(t,m),\xi}$ that
$t_1^{[\eta_l \leq]}$ is disjoint to, then by a probabilistic argument,
$\frac{1}{4^{l+1}}<(1-\frac{1}{2^{k+j+l}})^x$.

As $(1-\frac{1}{2^{k+j+l}})^{2k+j+l}<\frac{1}{e}<\frac{1}{2}$, it
follows that $x<2^{k+j+l}(2l+2)$, so we have at most $(2^{k+j+l}(2l+2))^2$
probematic pairs of $(\xi,\zeta)$ for a given pair of $(l,j)$. Therefore,
the number or problematic $(l,j,\xi,\zeta)$ is at most $\underset{l,j<k}{\Sigma}(2^{k+j+l}(2l+2))^2<2^{999k}$,
so by letting $\mu(k)=2^{2^{2^{999k}}}$ we're done. $\square$

\textbf{The following claim corresponds to step IV on page 6.}

\textbf{Claim 6. }There is a full system in $L[x_*]$.

\textbf{Proof: }The existence of such a system can be described by
a sentence $\psi$ in $\mathcal{L}_{\omega_1,\omega}$, and by Keisler's
completeness theorem it's absolute. By the previous claim, forcing
with finite systems over $L[x_*]$ preserves $\aleph_1$, hence we
can get a full system in $L[x_*]$.

\textbf{The following claim corresponds to step V on page 6.}

\textbf{Definition 7. }Fix a full system $S$. We define the formulas
$\phi_{rd}(x)$ and $\phi'_{rd}(x)$ (and similarly, $\phi_{gr}(x)$
and $\phi'_{gr}(x)$) as follows:

1. $\phi_{rd}(x)$ holds iff:

a. There is a tree $T_0$ which is a poor man generic tree over $L[x_*]$
(see see clause (A)(2) of step I in the above proof), so there is
$(n(k) : k<\omega)$ such that $(t(k),m(k))=(T_0 \restriction 2^{<n(k)}, ms_T \restriction 2^{<n(k)}) \in M_k$
(where for a closed tree $T$, the function $ms_T$ is defined as
$ms_T(\nu)=\mu(lim(T) \cap (2^{\omega})^{[\eta \leq]})$).

b. There is a tree $T_1$ which is a poor man generic tree over $L[x_*,T_0]$
(see clause (A)(1) of step I in the aboove proof), so $(t_n,m_n)=(T_1 \restriction 2^{\leq n}, ms_{T_1} \restriction 2^{\leq n}) \in N_n$
for every $n<\omega$.

c. For every $k<\omega$, $\rho_{rd}^S(t_{\lambda(t(k),m(k))},m_{\lambda(t(k),m(k))},t(k),m(k)) \leq x$.

2. $\phi'_{rd}(x)$ iff there is $y$ such that $\phi_{rd}(y)$ and
$x(n)=y(n)$ for $n$ large enough.

\textbf{Claim 8. }There above formulas are $\Sigma_3^1$.

\textbf{Proof: }Being contructible from $x_*$ is $\Sigma_2^1$, hence
{}``for every $G_{\delta}$ set $B$ of measure $0$, $B\cap T=\emptyset$
or $B$ is not contructible from $x_*$'' is $\Pi_2^1$ and the conclusion
follows. $\square$

\textbf{The following claim corresponds to step VI on page 6.}

\textbf{Claim 9. }$\phi'_{rd}(x)$ and $\phi'_{gr(x)}$ are contradictory.

\textbf{Proof: }Define a coloring of $[\omega_1]^2$ by $x(h(\xi,\zeta))$
for $\xi,\zeta<\omega_1$. If $\phi_{rd}(x)$, there are $T_0$, $T_1$
and $(n(k) : k<\omega)$ witnessing it. For $j<\omega$, $\eta_l \in T_1$
and $\alpha<\omega$ let $A_{j,l,n}$ be the set of $\xi<\omega_1$
such that $T_1^{\eta_l \leq}$ is disjoint to $\underset{l,k<\omega}{\cup}C_{l,j}^{(t(k),m(k)),\xi}$
and $f_{t_{\lambda(t(k),m(k))},m_{\lambda(t(k),m(k))}}^{(t(k),m(k)),\xi}(\eta_l,j)=\alpha$
for large enough $k$. This is a partition of $\omega_1$ to countably
many homogeneously red sets. Similarly, from $\phi_{gr}(x)$ we get
a partition of $\omega_1$ to countably many homogenously green sets,
so we get a contradiction.

Now suppose that $\phi_{rd}(x)$, $\phi_{gr}(y)$ and $x(n)=y(n)$
for $n>n^*$. There is a homogenously red set $A$ for $x$ and a
homogenously green set $B$ for $y$ such that $A\cap B$ is uncountable.
There is an infinite set $\{\xi_n : n<\omega\} \subseteq A\cap B$
such that $h(\xi_{n_1},\xi_{n_2})<h(\xi_{n_2},\xi_{n_3})$ has a fixed
truth value for $n_1<n_2<n_3$. By definition 6(A), $h(\xi_n,\xi_{n+1})$
is strictly increasing, hence it's $>n^*$ for $n$ large enough.
Therefore, for $n$ large enough, $red=x(h(\xi_n,\xi_{n+1}))=y(h(\xi_n,\xi_{n+1}))=green$,
a contradiction. $\square$

\textbf{Claim 13: }$A_{rd}=\{x : \phi'_{rd}(x)\}$ is not of measure
$0$.

\textbf{Proof: }Suppose that $b$ is a code for a $G_{\delta}$ set
of measure zero covering $A_{rd}$, then we get a poor man generic
tree $T_0$ over $L[x_*,b]$ and a poor man generic tree $T_1$ over
$L[x_*,b,T_0]$ (see step I in the above proof). Now let $x\in lim(T_0)$
such that $T_0$ and $T_1$ witness $\phi'_{rd}(x)$, then $x$ is
in no measure zero set coded in $L[x_*,b]$, contradicting the fact
that $x\in A_{rd}$ which is covered by the set coded by $b$. $\square$

\textbf{Claim 14: }$A_{rd}$ is not measurable. 

\textbf{Proof: }By the previous claim, its measure is not zero. By
the definition of $\phi'_{rd}$, the measure of $\{x : \phi'_{rd}(x),  \eta \leq x\}$
$(\eta \in 2^{<\omega})$ is determined by $lg(\eta)$. Therefore
$A_{rd}$ has outer measure $1$, and similarly for $A_{gr}$. As
they're disjoint, we get a contradiction. $\square$

\textbf{\large References}{\large \par}

{[}HwSh1067{]} Haim Horowitz and Saharon Shelah, Saccharinity with
ccc, preprint

{[}Sh176{]} Saharon Shelah, Can you take Solovay's inaccessible away?
Israel J. Math. 48 (1984), no. 1, 1-47

{[}So{]} Robert M. Solovay, A model of set theory in which every set
of reals is Lebesgue measurable, AM 92 (1970), 1-56

$\\$

(Haim Horowitz) Einstein Institute of Mathematics

Edmond J. Safra Campus,

The Hebrew University of Jerusalem.

Givat Ram, Jerusalem, 91904, Israel.

E-mail address: haim.horowitz@mail.huji.ac.il

$\\$

(Saharon Shelah) Einstein Institute of Mathematics

Edmond J. Safra Campus,

The Hebrew University of Jerusalem.

Givat Ram, Jerusalem, 91904, Israel.

Department of Mathematics

Hill Center - Busch Campus,

Rutgers, The State University of New Jersey.

110 Frelinghuysen road, Piscataway, NJ 08854-8019 USA

E-mail address: shelah@math.huji.ac.il
\end{document}